\documentclass[10pt]{amsart}
\usepackage{amsmath,amssymb,amsthm,url}
\usepackage[numbers]{natbib}
\usepackage{graphicx}

\newcommand{\rad}{\operatorname{rad}}
\newcommand{\Hopf}{\operatorname{Hopf}}

\theoremstyle{definition}

\newtheorem{theo}{Theorem}

\newtheorem{coro}{Corollary}

\begin{document}
\bibliographystyle{alpha}
\title[Extra structure and the WRT TQFT] {Extra structure and the Universal Construction for the Witten-Reshetikhin-Turaev TQFT}
 
\author[Patrick Gilmer]{Patrick M. Gilmer}
\address{Mathematics Department,
Louisiana State University,
Baton Rouge, Louisiana}
\email{gilmer@math.lsu.edu}
 \thanks{The first author was partially supported by  NSF-DMS-0905736}
\urladdr{\url{www.math.lsu.edu/~gilmer}}

\author{Xuanye Wang}
\address{Mathematics Department,
Louisiana State University,
Baton Rouge, Louisiana}
\email{xwang37@math.lsu.edu}

\date{\today}


\begin{abstract}
A TQFT is a functor from a cobordism category to the category of  vector spaces, satisfying certain properties. 
An important property is that the vector spaces
should be finite dimensional. For the WRT TQFT, the
relevant  $2+1$-cobordism category is 
built from manifolds which are equipped with 
an extra structure such as a $p_1$-structure, or an extended manifold
structure.
We perform the universal construction of \cite{BHMV2} on a cobordism category  without this extra structure and show that the resulting quantization functor
assigns an infinite dimensional vector space to the torus.
\end{abstract}

\maketitle

\section{Introduction}
A TQFT in dimension $2+1$ is a covariant functor $(V,Z)$ from 
some $(2+1)$-cobordism category $\mathcal{C}$ 
to the category of finite dimensional complex vector spaces  which assigns to the empty object  the vector space $\mathbb{C}$. 
Other properties  are usually required for a TQFT, and other ground rings are sometimes allowed. But for our purposes, this will do. 

Recall an object $\Sigma$ in $\mathcal{C}$ is a closed oriented surface with possibly some specified extra structure, and a morphism  
$C$ from $\Sigma_1$ to $\Sigma_2$ is an equivalence class of 
cobordisms from $\Sigma_1$ to $\Sigma_2$. Such a cobordism  can be loosely viewed as a compact oriented $3$-manifold (again possibly with some appropriate extra structure) with a boundary decomposed  
into  an incoming surface $-\Sigma_1$ and an  outgoing surface $\Sigma_2$.
Two cobordisms are considered equivalent if there is a orientation-preserving (extra structure preserving)  diffeomorphism between them which restricts to the identity on the boundary. 
Then $(V,Z)$ assigns a vector space $V(\Sigma)$ to
an object $\Sigma$, and a linear map $Z_{C}:V(\Sigma_1)\rightarrow V(\Sigma_2)$ to 
a  morphism $C:\Sigma_1 \rightarrow \Sigma_2$. 

The WRT-invariant is a $3$-manifold invariant which was first described by Witten in \cite{Witten} and  then rigorously defined by Reshetikhin and Turaev with quantum groups in \cite{RT}.  The approach to this invariant that we will use was developed by Blanchet, Habegger, Masbaum, and Vogel in \cite{BHMV1} 
with skein theory 
and then used by them to construct \cite{BHMV2} a TQFT on a 2+1 cobordism category where the objects and morphisms are equipped with $p_1$-structures. The question that we
consider in this paper is whether this construction 
based on the WRT-invariant still yields a TQFT when the extra structure is removed from the cobordism category.  The answer is no, as the resulting vector space associated to the torus has infinite dimension. See Theorem \ref{mainthm}. To be more precise: we follow this  construction after  assigning to each closed $3$-manifold the invariant of this 3-manifold equipped with a certain choice of extra structure. Our choice, which seems to us to be the most natural, is described in the next paragraph.

For each integer $p \ge 5$,   consider the complex valued invariant, $\langle \ \rangle_p$ of closed 
oriented $3$-manifolds equipped with a $p_1$-structure defined in \cite{BHMV2}. Here we must choose a primitive $2p$th root of unity $A\in \mathbb{C}$, and scalar $\kappa\in \mathbb{C}$ with $\kappa^6= A^{-6-p(p+1)/2}$.  
One may remove the dependence on this extra structure by defining  $\langle M\rangle'_p =  \langle \check M\rangle_p $ where $\check M$ is $M$ equipped with a $p_1$-structure with $\sigma$-invariant zero. See \cite[Appendix B]{BHMV2} for the definition of the $\sigma$-invariant.
If one uses extended manifold structures in lieu of $p_1$-structures as in \cite{Walker, Turaev, GM},  one would 
instead choose $\check M$ to have weight zero.

If $M$ is obtained by surgery to $S^3$ along a framed link $L$, then
\begin{equation}\label{invariant} 
\langle M\rangle'_p=\eta \mu^{-\sigma(L)} L(\omega_p). \notag
\end{equation}
Here we let  $\mu =\kappa^3$   and $\sigma(L)$ stands for the signature of the linking matrix of framed link $L$. 
Also $\omega_p$ is the skein specified in \cite[p.898]{BHMV2}, $\eta$ is the scalar as given in \cite[p.897]{BHMV2} and $L(\omega_p)$ is the Kauffman bracket of the cabling of $L$ by $\omega_p$.  One can easily extend this definition to the disconnected case by letting $\langle M_1 \sqcup M_2\rangle'_p=\langle M_1\rangle'_p \langle M_2\rangle'_p$.

A quantization functor is a covariant functor $(V,Z)$ from  $\mathcal{C}$ 
to a category of (not necessarily finite dimensional) complex vector spaces. 
Like a TQFT, it should assign to the empty object  the vector space $\mathbb{C}$. A certain naturally defined Hermitian form on $V(\Sigma)$ must also be non-degenerate.
One has that $\langle \ \rangle'_p$ is multiplicative and involutive. So we can perform the universal construction described in \cite[Prop. 1.1]{BHMV2} to construct a quantization functor  
from the ordinary $(2+1)$-cobordism category  (without any extra structure), which we will denote by 
$\mathcal{C}'$,
to the category of complex vector spaces

This is how the universal construction goes 
(when applied to $\mathcal{C}'$ and $\langle \ \rangle'_p$): 
 Given 
 an object $\Sigma$ in $\mathcal{C}'$, denote $\mathcal{V}'_{p}(\Sigma)$ as the vector space spanned by all compact oriented $3$-manifolds with boundary $\Sigma$ (or equivalently all  cobordisms $\{M:\emptyset \rightarrow \Sigma\}$). 
 There is  a hermitian form $\langle \ ,\ \rangle'_\Sigma$ on $\mathcal{V}'_p(\Sigma)$;
 This is specified on generators by
 \begin{equation}\label{UC}
\langle M,N\rangle'_\Sigma=\langle M \cup_\Sigma -N\rangle'_p  \notag
\end{equation}  
and then  extended sesquilinearly.  
Let $\rad \langle \ ,\ \rangle'_\Sigma$  denote the radical of the hermitian form $\langle \ , \ \rangle'_\Sigma$.
Define $V'_p(\Sigma)$ to be $\mathcal{V}'_p(\Sigma)\slash \rad \langle \ ,\ \rangle'_\Sigma$.  Given a morphism ${C}:\Sigma_1\rightarrow \Sigma_2$, define $\mathcal{Z}'_{p,{C}}: \mathcal{V}'_p(\Sigma_1)\rightarrow \mathcal{V}'_p(\Sigma_2)$ by assigning ${C} \cup_{\Sigma_1} {N}$ to any ${N} \in \mathcal{V}'_p(\Sigma_1)$ and extending linearly. Note that $\mathcal{Z}'_{p,{C}}$ send $\rad \langle \ ,\ \rangle'_{\Sigma_1}$ into $\rad \langle \ ,\ \rangle'_{\Sigma_2}$. So it induces a linear map $Z'_{p,{C}}: V'_p(\Sigma_1) \rightarrow V'_p(\Sigma_2).$  Then the quantization functor $(V'_p,Z'_p)$ is the rule  assigning $V'_p(\Sigma)$ to $\Sigma$ and $Z'_{p,{C}}$ to ${C}$.

Let $\mathcal{S}$ be a standard unknotted solid torus in 3-space, and let $T^2$ denote the boundary of  $\mathcal{S}$.
Let $w_i$ denote the $3$-manifold obtained by doing surgery to $\mathcal{S}$  along $i$ parallel copies of the core of $\mathcal{S}$ with framing $+1$. 
  Let $z_j$ denote the $3$-manifold obtained by doing surgery along 
the core of $\mathcal{S}$ with framing $j$.
We have the following theorem which will be proved in the  next section.

\begin{theo}\label{mainthm} 
For all $p \geq 5$, $V'_p(T^2)$ is infinite-dimensional. An infinite set of linearly independent elements in $V'_p(T^2)$ can be given by either $\{w_{2pk}\}$ or $\{z_{2pl}\}$, 
where $k$ and $l$
vary through the positive integers.      
Hence $V'_p$ is not a TQFT.
\end{theo}

If $(V,Z)$ is a quantization functor resulting from the universal construction,
then \cite[p.886]{BHMV2} there is a natural map 
$$t^{V}_{\Sigma_1, \Sigma_2}: V(\Sigma_1) \otimes V( \Sigma_2)  \rightarrow V(\Sigma_1 \sqcup \Sigma_2).$$ 
It is easy to see that this map must be injective.
Quinn \cite[Prop. 7.2]{Quinn} gave an argument that shows that the finite-dimensionality of $V(\Sigma)$ is implied by 
the functoriality of $V$ applied to 
a ``snake-shaped'' composition of cobordisms built from copies of $\Sigma \times I$ and the assumption that $t^{V}_{\Sigma, \Sigma}$ is an isomorphism. See also 
Kock \cite[Corollary 1.2.28]{Kock}.
This argument shows: 

\begin{coro}\label{maincoro}
The natural map $t^{{V'_p}}_{T^2, T^2}$ is not surjective. 
\end{coro}

\section{Proof of  Theorem}
 We first construct a bilinear form $B_p(\ ,\ )$ on  
$V'_p(T^2)$, then write out the $n\times n$ truncated matrix associated to $B_p(\ ,\ )$ with respect to $\{w_{2pk}\}$ and $\{z_{2pl}\}$ as $k,$ and  $l$ range from $1$ to $n$. 
 We will show there are infinitely many integers $n$ such that the truncated matrix of size $n \times n$ is nonsingular. Hence $\{w_{2pk}\}$ and $\{z_{2pl}\}$ are linearly independent and $V'_p(T^2)$ is infinite-dimensional.

\subsection{Bilinear form on $V'_p(T^2)$}

Every closed orientable connected $3$-manifold can be obtained by doing Dehn surgery in $S^3$, see \cite {lickorish}. This result can be used to show the following well-known related fact:
every connected orientable $3$-manifold with boundary $T^2$ can be obtained by doing surgery along some framed link $L$ in the solid torus 
$\mathcal{S}$. 
We denote the result of this surgery by 
$\mathcal{S}(L)$. 

According to the universal construction, elements in $V'_p(T^2)$ are represented by linear combinations of connected manifolds with boundary $T^2$.  
 Given two elements $w$ and $z$ in $V'_p(T^2)$ represented by $\mathcal{S}(L_w)$ and 
$\mathcal{S}(L_z)$, we glue together $\mathcal{S}(L_w)$ to  $\mathcal{S}(L_z)$  by the map on their boundary tori which switches the meridian and the longtitude (note this map is orientation reversing). 
The resulting  manifold is obtained by performing surgery on the 3-sphere along a framed link $(L_w,L_z)_{\Hopf}$ obtained by cabling the Hopf link with $L_w$ on one component  and $L_z$ on the other component. We define 
\begin{equation}
B_p(w,z)=\langle S^3((L_w,L_z)_{\Hopf})\rangle'_p \ .
\notag
\end{equation}
This extends to a well-defined symmetric bilinear form $B_p: V'_p(T^2) \times V'_p(T^2) \rightarrow \mathbb{C},$ as elements in $\rad \langle \ ,\ \rangle'_{T^2}$ will pair with any other element to give zero.

\subsection{Truncated square matrices}

\begin{figure}
\centering
\includegraphics[height=5.5cm]{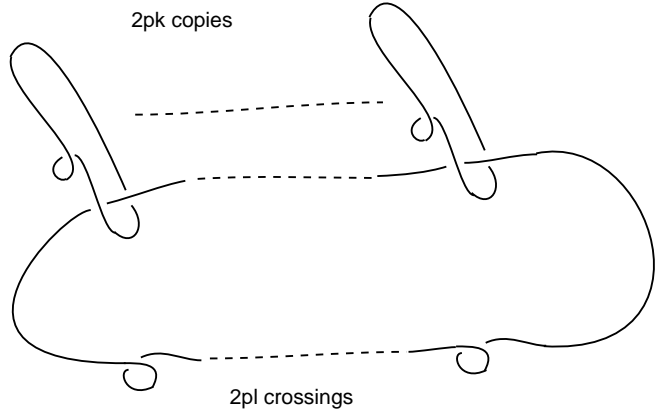}
\caption{Framed link  
$(L_{w_{2pk}},L_{z_{2pl}})_{\Hopf}$
}
\centering
\end{figure}
Note that 
$B_p(w_{2pk},z_{2pl})$ is Kauffman bracket of the  
$3$-manifold obtained by doing sugery to $S^3$ along the framed link  pictured in Figure 1. Blowing down \cite{Kirby} the $2pk$ unknotted components with framing $+1$, one by one,  and reducing the framing of the single component linked with these, we get an unknot $L'$
with framing $2p(l-k)$. Doing this surgery gives us the lens space $L(2p(l-k),1)$.

Thus $B(w_{2pk},z_{2pl})=\langle L(2p(l-k),1)\rangle'_p$.    
We let $s$ abbreviate  $l-k$ to simplify our expressions.   
To compute $\langle L(2ps,1)\rangle'_p$, we need to compute $U(t^{2ps} \omega_p)$. Here $t$ is the linear map from Kauffman skein module $K(S^1\times D^2)$ to itself induced by a positive twist 
and $U$ stands for an unknot with framing zero. 
According to \cite{BHMV1}, $t^{2ps} e_i= u_{k}^{2ps} e_i$, where $u_k=(-A)^{k(k+2)}$ and the $e_i$ are certain generators for $K(S^1\times D^2)$. 

As $A$ is a 2p-th root of unity, and  $u_k^{2ps}=1$, it follows that $t^{2ps}$ is the identity on $K(S^1\times D^2)$.    
We obtain
\begin{equation}
\langle L(2ps,1)\rangle'_p=\eta  U(t^{2ps} \omega_p) \mu^{-\sigma(L')} =\eta U(\omega_p) \mu^{-\sigma(L')}= \mu^{-\sigma(L')}
\ .
\notag
\end{equation}
This last equation holds as $\eta U(\omega_p) =1$ \cite[p.897]{BHMV2}.  Alternatively, 
$$\eta  U(\omega_p) = \langle S^1 \times S^2 \rangle'_p= \dim (V_p(S^2))=1.$$

Consider the matrix  with entries $B(w_{2pk},z_{2pl})$ as $k,l$ range over 
the integers from $1$ to $n$ (actually any  set of  $n$ integers  in increasing order will do). We will call this a truncated matrix of size $n$.
Then each entry in a truncated 
matrix just depends on the signature of the 
linking matrix of $L'$. Hence 
every truncated matrix has the form: 
\begin{equation}\label{matrix1}
\begin{bmatrix}
1 & \mu & \mu & \ldots & \mu\\
\mu^{-1} & 1 & \mu & \ldots & \mu\\
\mu^{-1} & \mu^{-1} & 1 & \ldots &\mu\\
\vdots & \vdots & \vdots & \ddots & \vdots\\
\mu^{-1} & \mu^{-1} & \mu^{-1} & \ldots & 1\\
\end{bmatrix}
\ .
\tag{*}
\end{equation}

\subsection{Determinants of the truncated matrices}
We will not try to show that the truncated 
  matrix of every 
 size
 is nonsingular. In fact, the truncated $2 \times 2$ matrix  has determinant zero. Instead, we show that two truncated  
 matrices of consecutive sizes can not be both singular. Therefore infinitely many nonsingular truncated matrices  
 exist. 

Consider the  $m\times m $ matrix $\mathcal{B}(a,m)$:   
\begin{equation}\label{matrix2}
\begin{bmatrix}
a & a-(1-\mu) & a-(1-\mu) & \ldots & a-(1-\mu)\\
\mu^{-1} & 1 & \mu & \ldots & \mu\\
\mu^{-1} & \mu^{-1} & 1 & \ldots &\mu\\
\vdots & \vdots & \vdots & \ddots & \vdots\\
\mu^{-1} & \mu^{-1} & \mu^{-1} & \ldots & 1\\
\end{bmatrix}
\ .
\notag
\end{equation}
Note that $\mu \neq 1$ for $p\ge 5$. We can apply the following elementary row and column operations to $\mathcal{B}(a,m)$: subtract the second column from the first column, add $-\mu^{-1}$ times the first row to the second row, then clearing all but the first entry in the first row, and obtain $(1-\mu) \oplus \mathcal{B}(f(a),m-1)$. Here $f(a)=\mu^{-1}(1-a)$. Applying same operations to the $\mathcal{B}(f(a),m-1)$ part of $(1-\mu) \oplus \mathcal{B}(f(a),m-1)$, we obtain $(1-\mu)I_2 \oplus \mathcal{B}(f^2(a),m-2)$. Repeating this $q$ times for some $q\leq m-1$, we see that $(1-\mu)I_q \oplus \mathcal{B}(f^q(a),m-q)$ is equivalent to $\mathcal{B}(a,m)$. Here we say two matrices are equivalent if they are related by a sequence of determinant preserving elementary row and column operations.

Note that $\mathcal{B}(1,n)$ is exactly the matrix \eqref{matrix1} of size $n$. Following from the above argument, it is clear that $\mathcal{B}(1,n)$ is equivalent to $(1-\mu)I_{n-1}\oplus \mathcal{B}(f^{n-1}(1),1)$ and that $\mathcal{B}(1,n+1)$ is equivalent to $(1-\mu)I_{n-1}\oplus \mathcal{B}(f^{n-1}(1),2)$. So both $\det \mathcal{B}(1,n)$ and $\det \mathcal{B}(1,n+1)$ can be written in terms of $f^{n-1}(1)$ as $\det \mathcal{B}(1,n)=(1-\mu)^{n-1}f^{n-1}(1)$ and $\det \mathcal{B}(1,n+1)=(1-\mu)^{n-1}(f^{n-1}(1)(1-\mu^{-1})+(\mu^{-1}-1))$. Therefore we have
\begin{equation}\label{det}
\det \mathcal{B}(1,n+1)=\det \mathcal{B}(1,n)(1-\mu^{-1})+(1-\mu)^{n-1}(\mu^{-1}-1)
\ .
\notag
\end{equation}

As a result, as long as $\mu \neq 1$, we can not have two consecutive singular truncated matrices since the term $(1-\mu)^{n-1}(\mu^{-1}-1)$ 
is non-zero.  
This completes the proof of Theorem \ref{mainthm}.

\bibliography{Version6}

\end{document}